\documentclass{gtart}
\def\ifplaintex{\expandafter\ifx\csname documentclass\endcsname\relax}

\def\gtp{{\mathsurround=0pt\it $\cal G\mskip-2mu$eometry \&\ 
$\cal T\!\!$opology $\cal P\!$ublications}}  % GT publications

\def\recd{{\small Received:\qua\receiveddate\ifx\reviseddate\relax
\else\qquad Revised:\qua\reviseddate\fi\par}} 

\def\lognumber#1{\def\thelognumber{#1}}
\def\volumenumber#1{\def\thevolumenumber{#1}}
\def\volumeyear#1{\def\thevolumeyear{#1}}
\def\papernumber#1{\def\thepapernumber{#1}}
\def\pagenumbers#1#2{\def\startpage{#1}\def\finishpage{#2}}
\def\published#1{\def\publishdate{#1}}

\def\received#1{\def\receiveddate{#1}}

\def\accepted#1{\def\accepteddate{#1}}

\long\def\asciiabstract#1{\long\def\theasciiabstract{#1}}

\let\\\par\let\thelognumber\relax\let\thevolumenumber\relax
\let\thepapernumber\relax\let\thevolumeyear\relax\let\startpage\relax
\let\finishpage\relax\let\publishdate\relax\let\receiveddate\relax
\let\reviseddate\relax\let\accepteddate\relax\let\theasciititle\relax
\let\theasciiauthors\relax
\let\theasciiabstract\relax

\let\theasciiemail\relax

\ifplaintex
\font\logobig=cmssbx10 scaled 3836
\font\logomed=cmssbx10 scaled 2557
\else
\font\logobig=cmssbx10 scaled 4200
\font\logomed=cmssbx10 scaled 2800
\fi

\long\def\makeagttitle{   %%% start of definition of \makeagttitle
\count0=\startpage
\agt\hfill      %   Journal title (top left) 
\hbox to 45truept{\vbox to 0pt{\vglue -13truept{\logomed A\kern -.37em{\logobig 
T}\kern -.38em G}\vss}\hss}
\break
{\small Volume \thevolumenumber\ (\thevolumeyear)
\startpage--\finishpage\nl
Published: \publishdate}

\vglue .25truein

{\parskip=0pt\leftskip 0pt plus
1fil\def\\{\par\smallskip}{\Large\bf\thetitle}\par\medskip} \vglue
0.05truein

{\parskip=0pt\leftskip 0pt plus 1fil\def\\{\par}{\sc\theauthors}
\par\medskip}%
 
\vglue 0.03truein 

{\small\leftskip 25truept\rightskip 25truept{\bf Abstract}\stdspace\theabstract

{\bf AMS Classification}\stdspace\theprimaryclass
\ifx\thesecondaryclass\relax\else; \thesecondaryclass\fi\par
{\bf Keywords}\stdspace \thekeywords\par}\vglue 7truept

}   %%%% end of definition of \makeagttitle

\ifplaintex
\font\phead=cmsl9 scaled 950
\font\pnum=cmbx10 scaled 913
\font\pfoot=cmsl9 scaled 950
\headline{\vbox to 0pt{\vskip -4.5mm\line{\small\phead\ifnum
\count0=\startpage ISSN 1472-2739 (on-line) 1472-2747 (printed)
\hfill {\pnum\folio}\else\ifodd\count0\def\\{ }% 
\ifx\theshorttitle\relax\thetitle\else\theshorttitle\fi\hfill{\pnum\folio}
\else\def\\{ and }{\pnum\folio}\hfill\ifx\theshortauthors\relax\theauthors
\else\theshortauthors\fi\fi\fi}\vss}}
\footline{\vbox to 0pt{\vglue 0mm\line{\small\pfoot\ifnum\count0=\startpage
\copyright\ \gtp\hfill\else
\agt, Volume \thevolumenumber\ (\thevolumeyear)\hfill\fi}\vss}}
\else
\font=cmsl9 scaled 1050
\font\lnum=cmbx10 
\font=cmsl9 scaled 1050
\makeatletter
\def\@oddhead{{\small\lhead\ifnum\count0=\startpage ISSN 1472-2739 
(on-line) 1472-2747 (printed)\hfill {\lnum\number\count0}\else\ifodd\count0
\def\\{ }\ifx\theshorttitle\relax \thetitle \else\theshorttitle\fi\hfill
{\lnum\number\count0}\else\def\\{ and }{\lnum\number\count0}
\hfill\ifx\theshortauthors\relax 
\theauthors\else\theshortauthors\fi\fi\fi}}\def\@evenhead{\@oddhead}
\def\@oddfoot{\small\lfoot\ifnum\count0=\startpage\copyright\ \gtp\hfill\else
\agt, Volume \thevolumenumber\ (\thevolumeyear)\hfill\fi}
\def\@evenfoot{\@oddfoot}
\makeatother
\fi
\let\maketitlepage\makeagttitle
\let\makeshorttitle\maketitlepage
\let\maketitle\maketitlepage

\newwrite\gtoutfile
\long\gdef\makeheadfile{  %%% start of definition of \makeheadfile
{\def\\{, }\def\s{ }
\immediate\openout\gtoutfile head.xxx
\immediate\write\gtoutfile{To: math@arxiv.org}
\immediate\write\gtoutfile{Subject: put OR rep NNNNN:ppppp}
\immediate\write\gtoutfile{--text follows this line--}
\immediate\write\gtoutfile{Proxy-for: \ifx\theasciiauthors\relax
\theauthors\else\theasciiauthors\fi\s<\ifx\theasciiemail\relax\theemail\else\theasciiemail\fi>}
\immediate\write\gtoutfile{\noexpand\\}
\immediate\write\gtoutfile{Authors: \ifx\theasciiauthors\relax
\theauthors\else\theasciiauthors\fi}
{\def\\{ }\immediate\write\gtoutfile{Title: \ifx\theasciititle\relax
\thetitle\else\theasciititle\fi}}
\immediate\write\gtoutfile{Subj-class: GT or SG, GR etc}
\immediate\write\gtoutfile{MSC-class: \theprimaryclass\ifx\thesecondaryclass\relax\else, \thesecondaryclass\fi}
\immediate\write\gtoutfile{Journal-ref: Algebr. Geom. Topol. \thevolumenumber\s
(\thevolumeyear) \startpage-\finishpage}
\immediate\write\gtoutfile{Comments: Published by Algebraic and
Geometric Topology at}
\immediate\write\gtoutfile{\s\s\s  http://www.maths.warwick.ac.uk/agt/AGTVol\thevolumenumber/agt-\thevolumenumber-\thepapernumber.abs.html}
\immediate\write\gtoutfile{\noexpand\\}
\immediate\write\gtoutfile{}
\ifx\theasciiabstract\relax
\immediate\write\gtoutfile{\theabstract}\else
\immediate\write\gtoutfile{\theasciiabstract}\fi
\immediate\write\gtoutfile{}
\immediate\write\gtoutfile{\noexpand\\}
\immediate\write\gtoutfile{}
\immediate\closeout\gtoutfile}}  %%% end of definition of \makeheadfile

\def\maketitlepage{\makeagttitle\makeheadfile}
\let\makeshorttitle\maketitlepage
\let\maketitle\maketitlepage

\def\ifplaintex{\expandafter\ifx\csname documentclass\endcsname\relax}

\def\gtp{{\mathsurround=0pt\it $\cal G\mskip-2mu$eometry \&\ 
$\cal T\!\!$opology $\cal P\!$ublications}}  % GT publications

\def\recd{{\small Received:\qua\receiveddate\ifx\reviseddate\relax
\else\qquad Revised:\qua\reviseddate\fi\par}} 

\def\lognumber#1{\def\thelognumber{#1}}
\def\volumenumber#1{\def\thevolumenumber{#1}}
\def\volumeyear#1{\def\thevolumeyear{#1}}
\def\papernumber#1{\def\thepapernumber{#1}}
\def\pagenumbers#1#2{\def\startpage{#1}\def\finishpage{#2}}
\def\published#1{\def\publishdate{#1}}

\def\received#1{\def\receiveddate{#1}}

\def\accepted#1{\def\accepteddate{#1}}

\long\def\asciiabstract#1{\long\def\theasciiabstract{#1}}

\let\\\par\let\thelognumber\relax\let\thevolumenumber\relax
\let\thepapernumber\relax\let\thevolumeyear\relax\let\startpage\relax
\let\finishpage\relax\let\publishdate\relax\let\receiveddate\relax
\let\reviseddate\relax\let\accepteddate\relax\let\theasciititle\relax
\let\theasciiauthors\relax
\let\theasciiabstract\relax

\let\theasciiemail\relax

\ifplaintex
\font\logobig=cmssbx10 scaled 3836
\font\logomed=cmssbx10 scaled 2557
\else
\font\logobig=cmssbx10 scaled 4200
\font\logomed=cmssbx10 scaled 2800
\fi

\long\def\makeagttitle{   %%% start of definition of \makeagttitle
\count0=\startpage
\agt\hfill      %   Journal title (top left) 
\hbox to 45truept{\vbox to 0pt{\vglue -13truept{\logomed A\kern -.37em{\logobig 
T}\kern -.38em G}\vss}\hss}
\break
{\small Volume \thevolumenumber\ (\thevolumeyear)
\startpage--\finishpage\nl
Published: \publishdate}

\vglue .25truein

{\parskip=0pt\leftskip 0pt plus
1fil\def\\{\par\smallskip}{\Large\bf\thetitle}\par\medskip} \vglue
0.05truein

{\parskip=0pt\leftskip 0pt plus 1fil\def\\{\par}{\sc\theauthors}
\par\medskip}%
 
\vglue 0.03truein 

{\small\leftskip 25truept\rightskip 25truept{\bf Abstract}\stdspace\theabstract

{\bf AMS Classification}\stdspace\theprimaryclass
\ifx\thesecondaryclass\relax\else; \thesecondaryclass\fi\par
{\bf Keywords}\stdspace \thekeywords\par}\vglue 7truept

}   %%%% end of definition of \makeagttitle

\ifplaintex
\font\phead=cmsl9 scaled 950
\font\pnum=cmbx10 scaled 913
\font\pfoot=cmsl9 scaled 950
\headline{\vbox to 0pt{\vskip -4.5mm\line{\small\phead\ifnum
\count0=\startpage ISSN 1472-2739 (on-line) 1472-2747 (printed)
\hfill {\pnum\folio}\else\ifodd\count0\def\\{ }% 
\ifx\theshorttitle\relax\thetitle\else\theshorttitle\fi\hfill{\pnum\folio}
\else\def\\{ and }{\pnum\folio}\hfill\ifx\theshortauthors\relax\theauthors
\else\theshortauthors\fi\fi\fi}\vss}}
\footline{\vbox to 0pt{\vglue 0mm\line{\small\pfoot\ifnum\count0=\startpage
\copyright\ \gtp\hfill\else
\agt, Volume \thevolumenumber\ (\thevolumeyear)\hfill\fi}\vss}}
\else
\font=cmsl9 scaled 1050
\font\lnum=cmbx10 
\font=cmsl9 scaled 1050
\makeatletter
\def\@oddhead{{\small\lhead\ifnum\count0=\startpage ISSN 1472-2739 
(on-line) 1472-2747 (printed)\hfill {\lnum\number\count0}\else\ifodd\count0
\def\\{ }\ifx\theshorttitle\relax \thetitle \else\theshorttitle\fi\hfill
{\lnum\number\count0}\else\def\\{ and }{\lnum\number\count0}
\hfill\ifx\theshortauthors\relax 
\theauthors\else\theshortauthors\fi\fi\fi}}\def\@evenhead{\@oddhead}
\def\@oddfoot{\small\lfoot\ifnum\count0=\startpage\copyright\ \gtp\hfill\else
\agt, Volume \thevolumenumber\ (\thevolumeyear)\hfill\fi}
\def\@evenfoot{\@oddfoot}
\makeatother
\fi
\let\maketitlepage\makeagttitle
\let\makeshorttitle\maketitlepage
\let\maketitle\maketitlepage

\newwrite\gtoutfile
\long\gdef\makeheadfile{  %%% start of definition of \makeheadfile
{\def\\{, }\def\s{ }
\immediate\openout\gtoutfile head.xxx
\immediate\write\gtoutfile{To: math@arxiv.org}
\immediate\write\gtoutfile{Subject: put OR rep NNNNN:ppppp}
\immediate\write\gtoutfile{--text follows this line--}
\immediate\write\gtoutfile{Proxy-for: \ifx\theasciiauthors\relax
\theauthors\else\theasciiauthors\fi\s<\ifx\theasciiemail\relax\theemail\else\theasciiemail\fi>}
\immediate\write\gtoutfile{\noexpand\\}
\immediate\write\gtoutfile{Authors: \ifx\theasciiauthors\relax
\theauthors\else\theasciiauthors\fi}
{\def\\{ }\immediate\write\gtoutfile{Title: \ifx\theasciititle\relax
\thetitle\else\theasciititle\fi}}
\immediate\write\gtoutfile{Subj-class: GT or SG, GR etc}
\immediate\write\gtoutfile{MSC-class: \theprimaryclass\ifx\thesecondaryclass\relax\else, \thesecondaryclass\fi}
\immediate\write\gtoutfile{Journal-ref: Algebr. Geom. Topol. \thevolumenumber\s
(\thevolumeyear) \startpage-\finishpage}
\immediate\write\gtoutfile{Comments: Published by Algebraic and
Geometric Topology at}
\immediate\write\gtoutfile{\s\s\s  http://www.maths.warwick.ac.uk/agt/AGTVol\thevolumenumber/agt-\thevolumenumber-\thepapernumber.abs.html}
\immediate\write\gtoutfile{\noexpand\\}
\immediate\write\gtoutfile{}
\ifx\theasciiabstract\relax
\immediate\write\gtoutfile{\theabstract}\else
\immediate\write\gtoutfile{\theasciiabstract}\fi
\immediate\write\gtoutfile{}
\immediate\write\gtoutfile{\noexpand\\}
\immediate\write\gtoutfile{}
\immediate\closeout\gtoutfile}}  %%% end of definition of \makeheadfile

\def\maketitlepage{\makeagttitle\makeheadfile}
\let\makeshorttitle\maketitlepage
\let\maketitle\maketitlepage

\def\ifplaintex{\expandafter\ifx\csname documentclass\endcsname\relax}

\def\gtp{{\mathsurround=0pt\it $\cal G\mskip-2mu$eometry \&\ 
$\cal T\!\!$opology $\cal P\!$ublications}}  % GT publications

\def\recd{{\small Received:\qua\receiveddate\ifx\reviseddate\relax
\else\qquad Revised:\qua\reviseddate\fi\par}} 

\def\lognumber#1{\def\thelognumber{#1}}
\def\volumenumber#1{\def\thevolumenumber{#1}}
\def\volumeyear#1{\def\thevolumeyear{#1}}
\def\papernumber#1{\def\thepapernumber{#1}}
\def\pagenumbers#1#2{\def\startpage{#1}\def\finishpage{#2}}
\def\published#1{\def\publishdate{#1}}

\def\received#1{\def\receiveddate{#1}}

\def\accepted#1{\def\accepteddate{#1}}

\long\def\asciiabstract#1{\long\def\theasciiabstract{#1}}

\let\\\par\let\thelognumber\relax\let\thevolumenumber\relax
\let\thepapernumber\relax\let\thevolumeyear\relax\let\startpage\relax
\let\finishpage\relax\let\publishdate\relax\let\receiveddate\relax
\let\reviseddate\relax\let\accepteddate\relax\let\theasciititle\relax
\let\theasciiauthors\relax
\let\theasciiabstract\relax

\let\theasciiemail\relax

\ifplaintex
\font\logobig=cmssbx10 scaled 3836
\font\logomed=cmssbx10 scaled 2557
\else
\font\logobig=cmssbx10 scaled 4200
\font\logomed=cmssbx10 scaled 2800
\fi

\long\def\makeagttitle{   %%% start of definition of \makeagttitle
\count0=\startpage
\agt\hfill      %   Journal title (top left) 
\hbox to 45truept{\vbox to 0pt{\vglue -13truept{\logomed A\kern -.37em{\logobig 
T}\kern -.38em G}\vss}\hss}
\break
{\small Volume \thevolumenumber\ (\thevolumeyear)
\startpage--\finishpage\nl
Published: \publishdate}

\vglue .25truein

{\parskip=0pt\leftskip 0pt plus
1fil\def\\{\par\smallskip}{\Large\bf\thetitle}\par\medskip} \vglue
0.05truein

{\parskip=0pt\leftskip 0pt plus 1fil\def\\{\par}{\sc\theauthors}
\par\medskip}%
 
\vglue 0.03truein 

{\small\leftskip 25truept\rightskip 25truept{\bf Abstract}\stdspace\theabstract

{\bf AMS Classification}\stdspace\theprimaryclass
\ifx\thesecondaryclass\relax\else; \thesecondaryclass\fi\par
{\bf Keywords}\stdspace \thekeywords\par}\vglue 7truept

}   %%%% end of definition of \makeagttitle

\ifplaintex
\font\phead=cmsl9 scaled 950
\font\pnum=cmbx10 scaled 913
\font\pfoot=cmsl9 scaled 950
\headline{\vbox to 0pt{\vskip -4.5mm\line{\small\phead\ifnum
\count0=\startpage ISSN 1472-2739 (on-line) 1472-2747 (printed)
\hfill {\pnum\folio}\else\ifodd\count0\def\\{ }% 
\ifx\theshorttitle\relax\thetitle\else\theshorttitle\fi\hfill{\pnum\folio}
\else\def\\{ and }{\pnum\folio}\hfill\ifx\theshortauthors\relax\theauthors
\else\theshortauthors\fi\fi\fi}\vss}}
\footline{\vbox to 0pt{\vglue 0mm\line{\small\pfoot\ifnum\count0=\startpage
\copyright\ \gtp\hfill\else
\agt, Volume \thevolumenumber\ (\thevolumeyear)\hfill\fi}\vss}}
\else
\font=cmsl9 scaled 1050
\font\lnum=cmbx10 
\font=cmsl9 scaled 1050
\makeatletter
\def\@oddhead{{\small\lhead\ifnum\count0=\startpage ISSN 1472-2739 
(on-line) 1472-2747 (printed)\hfill {\lnum\number\count0}\else\ifodd\count0
\def\\{ }\ifx\theshorttitle\relax \thetitle \else\theshorttitle\fi\hfill
{\lnum\number\count0}\else\def\\{ and }{\lnum\number\count0}
\hfill\ifx\theshortauthors\relax 
\theauthors\else\theshortauthors\fi\fi\fi}}\def\@evenhead{\@oddhead}
\def\@oddfoot{\small\lfoot\ifnum\count0=\startpage\copyright\ \gtp\hfill\else
\agt, Volume \thevolumenumber\ (\thevolumeyear)\hfill\fi}
\def\@evenfoot{\@oddfoot}
\makeatother
\fi
\let\maketitlepage\makeagttitle
\let\makeshorttitle\maketitlepage
\let\maketitle\maketitlepage

\newwrite\gtoutfile
\long\gdef\makeheadfile{  %%% start of definition of \makeheadfile
{\def\\{, }\def\s{ }
\immediate\openout\gtoutfile head.xxx
\immediate\write\gtoutfile{To: math@arxiv.org}
\immediate\write\gtoutfile{Subject: put OR rep NNNNN:ppppp}
\immediate\write\gtoutfile{--text follows this line--}
\immediate\write\gtoutfile{Proxy-for: \ifx\theasciiauthors\relax
\theauthors\else\theasciiauthors\fi\s<\ifx\theasciiemail\relax\theemail\else\theasciiemail\fi>}
\immediate\write\gtoutfile{\noexpand\\}
\immediate\write\gtoutfile{Authors: \ifx\theasciiauthors\relax
\theauthors\else\theasciiauthors\fi}
{\def\\{ }\immediate\write\gtoutfile{Title: \ifx\theasciititle\relax
\thetitle\else\theasciititle\fi}}
\immediate\write\gtoutfile{Subj-class: GT or SG, GR etc}
\immediate\write\gtoutfile{MSC-class: \theprimaryclass\ifx\thesecondaryclass\relax\else, \thesecondaryclass\fi}
\immediate\write\gtoutfile{Journal-ref: Algebr. Geom. Topol. \thevolumenumber\s
(\thevolumeyear) \startpage-\finishpage}
\immediate\write\gtoutfile{Comments: Published by Algebraic and
Geometric Topology at}
\immediate\write\gtoutfile{\s\s\s  http://www.maths.warwick.ac.uk/agt/AGTVol\thevolumenumber/agt-\thevolumenumber-\thepapernumber.abs.html}
\immediate\write\gtoutfile{\noexpand\\}
\immediate\write\gtoutfile{}
\ifx\theasciiabstract\relax
\immediate\write\gtoutfile{\theabstract}\else
\immediate\write\gtoutfile{\theasciiabstract}\fi
\immediate\write\gtoutfile{}
\immediate\write\gtoutfile{\noexpand\\}
\immediate\write\gtoutfile{}
\immediate\closeout\gtoutfile}}  %%% end of definition of \makeheadfile

\def\maketitlepage{\makeagttitle\makeheadfile}
\let\makeshorttitle\maketitlepage
\let\maketitle\maketitlepage

\def\ifplaintex{\expandafter\ifx\csname documentclass\endcsname\relax}

\def\gtp{{\mathsurround=0pt\it $\cal G\mskip-2mu$eometry \&\ 
$\cal T\!\!$opology $\cal P\!$ublications}}  % GT publications

\def\recd{{\small Received:\qua\receiveddate\ifx\reviseddate\relax
\else\qquad Revised:\qua\reviseddate\fi\par}} 

\def\lognumber#1{\def\thelognumber{#1}}
\def\volumenumber#1{\def\thevolumenumber{#1}}
\def\volumeyear#1{\def\thevolumeyear{#1}}
\def\papernumber#1{\def\thepapernumber{#1}}
\def\pagenumbers#1#2{\def\startpage{#1}\def\finishpage{#2}}
\def\published#1{\def\publishdate{#1}}

\def\received#1{\def\receiveddate{#1}}

\def\accepted#1{\def\accepteddate{#1}}

\long\def\asciiabstract#1{\long\def\theasciiabstract{#1}}

\let\\\par\let\thelognumber\relax\let\thevolumenumber\relax
\let\thepapernumber\relax\let\thevolumeyear\relax\let\startpage\relax
\let\finishpage\relax\let\publishdate\relax\let\receiveddate\relax
\let\reviseddate\relax\let\accepteddate\relax\let\theasciititle\relax
\let\theasciiauthors\relax
\let\theasciiabstract\relax

\let\theasciiemail\relax

\ifplaintex
\font\logobig=cmssbx10 scaled 3836
\font\logomed=cmssbx10 scaled 2557
\else
\font\logobig=cmssbx10 scaled 4200
\font\logomed=cmssbx10 scaled 2800
\fi

\long\def\makeagttitle{   %%% start of definition of \makeagttitle
\count0=\startpage
\agt\hfill      %   Journal title (top left) 
\hbox to 45truept{\vbox to 0pt{\vglue -13truept{\logomed A\kern -.37em{\logobig 
T}\kern -.38em G}\vss}\hss}
\break
{\small Volume \thevolumenumber\ (\thevolumeyear)
\startpage--\finishpage\nl
Published: \publishdate}

\vglue .25truein

{\parskip=0pt\leftskip 0pt plus
1fil\def\\{\par\smallskip}{\Large\bf\thetitle}\par\medskip} \vglue
0.05truein

{\parskip=0pt\leftskip 0pt plus 1fil\def\\{\par}{\sc\theauthors}
\par\medskip}%
 
\vglue 0.03truein 

{\small\leftskip 25truept\rightskip 25truept{\bf Abstract}\stdspace\theabstract

{\bf AMS Classification}\stdspace\theprimaryclass
\ifx\thesecondaryclass\relax\else; \thesecondaryclass\fi\par
{\bf Keywords}\stdspace \thekeywords\par}\vglue 7truept

}   %%%% end of definition of \makeagttitle

\ifplaintex
\font\phead=cmsl9 scaled 950
\font\pnum=cmbx10 scaled 913
\font\pfoot=cmsl9 scaled 950
\headline{\vbox to 0pt{\vskip -4.5mm\line{\small\phead\ifnum
\count0=\startpage ISSN 1472-2739 (on-line) 1472-2747 (printed)
\hfill {\pnum\folio}\else\ifodd\count0\def\\{ }% 
\ifx\theshorttitle\relax\thetitle\else\theshorttitle\fi\hfill{\pnum\folio}
\else\def\\{ and }{\pnum\folio}\hfill\ifx\theshortauthors\relax\theauthors
\else\theshortauthors\fi\fi\fi}\vss}}
\footline{\vbox to 0pt{\vglue 0mm\line{\small\pfoot\ifnum\count0=\startpage
\copyright\ \gtp\hfill\else
\agt, Volume \thevolumenumber\ (\thevolumeyear)\hfill\fi}\vss}}
\else
\font=cmsl9 scaled 1050
\font\lnum=cmbx10 
\font=cmsl9 scaled 1050
\makeatletter
\def\@oddhead{{\small\lhead\ifnum\count0=\startpage ISSN 1472-2739 
(on-line) 1472-2747 (printed)\hfill {\lnum\number\count0}\else\ifodd\count0
\def\\{ }\ifx\theshorttitle\relax \thetitle \else\theshorttitle\fi\hfill
{\lnum\number\count0}\else\def\\{ and }{\lnum\number\count0}
\hfill\ifx\theshortauthors\relax 
\theauthors\else\theshortauthors\fi\fi\fi}}\def\@evenhead{\@oddhead}
\def\@oddfoot{\small\lfoot\ifnum\count0=\startpage\copyright\ \gtp\hfill\else
\agt, Volume \thevolumenumber\ (\thevolumeyear)\hfill\fi}
\def\@evenfoot{\@oddfoot}
\makeatother
\fi
\let\maketitlepage\makeagttitle
\let\makeshorttitle\maketitlepage
\let\maketitle\maketitlepage

\newwrite\gtoutfile
\long\gdef\makeheadfile{  %%% start of definition of \makeheadfile
{\def\\{, }\def\s{ }
\immediate\openout\gtoutfile head.xxx
\immediate\write\gtoutfile{To: math@arxiv.org}
\immediate\write\gtoutfile{Subject: put OR rep NNNNN:ppppp}
\immediate\write\gtoutfile{--text follows this line--}
\immediate\write\gtoutfile{Proxy-for: \ifx\theasciiauthors\relax
\theauthors\else\theasciiauthors\fi\s<\ifx\theasciiemail\relax\theemail\else\theasciiemail\fi>}
\immediate\write\gtoutfile{\noexpand\\}
\immediate\write\gtoutfile{Authors: \ifx\theasciiauthors\relax
\theauthors\else\theasciiauthors\fi}
{\def\\{ }\immediate\write\gtoutfile{Title: \ifx\theasciititle\relax
\thetitle\else\theasciititle\fi}}
\immediate\write\gtoutfile{Subj-class: GT or SG, GR etc}
\immediate\write\gtoutfile{MSC-class: \theprimaryclass\ifx\thesecondaryclass\relax\else, \thesecondaryclass\fi}
\immediate\write\gtoutfile{Journal-ref: Algebr. Geom. Topol. \thevolumenumber\s
(\thevolumeyear) \startpage-\finishpage}
\immediate\write\gtoutfile{Comments: Published by Algebraic and
Geometric Topology at}
\immediate\write\gtoutfile{\s\s\s  http://www.maths.warwick.ac.uk/agt/AGTVol\thevolumenumber/agt-\thevolumenumber-\thepapernumber.abs.html}
\immediate\write\gtoutfile{\noexpand\\}
\immediate\write\gtoutfile{}
\ifx\theasciiabstract\relax
\immediate\write\gtoutfile{\theabstract}\else
\immediate\write\gtoutfile{\theasciiabstract}\fi
\immediate\write\gtoutfile{}
\immediate\write\gtoutfile{\noexpand\\}
\immediate\write\gtoutfile{}
\immediate\closeout\gtoutfile}}  %%% end of definition of \makeheadfile

\def\maketitlepage{\makeagttitle\makeheadfile}
\let\makeshorttitle\maketitlepage
\let\maketitle\maketitlepage

\lognumber{38}
\volumenumber{2}
\volumeyear{2002}
\papernumber{38}
\pagenumbers{937}{947}
\received{23 July 2002}
\accepted{3 October 2002}
\published{22 October 2002}

\usepackage{amssymb} 
\usepackage{amsmath}  %%% newams

\input xypic
\input xy
\xyoption{all}

\begin{document}

\newcommand{\colim}{\ensuremath{\mathrm{colim}}}
\newcommand{\ilim}{\ensuremath{\mathrm{lim}}}
\newcommand{\hocolim}{\ensuremath{\mathrm{hocolim}}}
\newcommand{\holim}{\ensuremath{\mathrm{holim}}}
\newtheorem{theor}{Theorem}[section]  
\newtheorem{lemma}[theor]{Lemma}     
\newtheorem{prop}[theor]{Proposition}
\newtheorem{coro}[theor]{Corollary}
\newtheorem{conj}[theor]{Conjecture}
\theoremstyle{definition}
\newtheorem{defi}[theor]{Definition}   
\numberwithin{equation}{section}

\title{On Real-oriented Johnson-Wilson cohomology}
\author{Po Hu}
\address{Department of Mathematics, Wayne State University\\656 
W. Kirby, Detroit, MI 48202, USA}
\email{po@math.wayne.edu}
\begin{abstract}
Answering a question of W. S. Wilson, I introduce a 
${\mathbb Z}/2$-equivariant Atiyah-Real analogue of Johnson-Wilson 
cohomology theory $BP\langle n\rangle$, whose coefficient ring is the 
$\leq n$-chromatic 
part of Landweber's Real cobordism ring.
\end{abstract}

\asciiabstract{Answering a question of W. S. Wilson, I introduce a
Z/2-equivariant Atiyah-Real analogue of Johnson-Wilson cohomology
theory BP<n>, whose coefficient ring is the =< n-chromatic part of
Landweber's Real cobordism ring.}


\maketitle

\section{Introduction}

Recall Johnson-Wilson's spectrum $BP \langle n\rangle$, 
constructed in~\cite{wilson}. The complex cobordism spectrum $MU$, localized 
at $2$, splits as a wedge-sum of suspensions of the Brown-Peterson 
spectrum $BP$~\cite{land}. We have $BP_{\ast} = {\mathbb Z}_{(2)}[v_1, 
v_2, \ldots]$, with $dim(v_i) = 2(2^i-1)$. For each $n$, the Johnson-Wilson 
spectrum $BP\langle n\rangle$ comes with a map $BP \rightarrow BP \langle n
\rangle$, and one has that 
\begin{equation}
BP \langle n\rangle_{\ast} = {\mathbb Z}_{(2)}[v_1, v_2, \ldots, v_n] .
\label{bpncoeff}
\end{equation}
In particular, $BP \langle n\rangle_{\ast}$ is a quotient ring of $BP_{\ast}$.
The fact that such $BP\langle n\rangle$ exists is, of course, 
today no longer surprising. In fact, one can construct 
$BP$ almost formally by ``killing a suitable regular ideal $I_n$ in 
$MU_{(2)}$'' (see~\cite{ekmm}).

In connection with certain questions on Lie groups (which will not be 
discussed here), Steve Wilson recently asked if the spectrum $BP\langle n
\rangle$ has a Landweber-Real analogue, i.e., if there exists a spectrum 
$BP{\mathbb R}\langle n\rangle$ whose coefficient ring is the 
quotient of Landweber's cobordism ring $M{\mathbb R}_{\star}$~(\cite{araki2}, \cite{araki3},
\cite{hk}, \cite{land}) by all elements ``not related to $v_0, \ldots, 
v_n$''. (Here, $M{\mathbb R}_{\star}$ denotes the $RO({\mathbb Z}/2)$-graded
coefficient ring, as opposed to the integer-graded coefficient 
ring.) This can be given an exact meaning, which I shall explain in the 
next section.
First, however, I shall describe, in general terms, the main result of this 
paper, and its contribution to the present state of the subject.

In this paper, I completely answer Steve Wilson's question in the 
affirmative. The construction of the spectrum $BP\mathbb{R} \langle n
\rangle$ is straightforward:
analogously to $MU$-theory, general tools are now avaliable in 
$M{\mathbb R}$-theory. In particular, there is an embedding 
$MU_{\ast} \rightarrow M{\mathbb R}_{\star}$ (in an appropriate sense), 
and it is possible to ``quotient out'' $M{\mathbb R}$ by an 
ideal of $MU_{\ast}$ using the tools of~\cite{ekmm}. This method is 
described in detail in~\cite{hk}. The construction of my spectrum 
$BP{\mathbb R}\langle n\rangle$ is that one 
simply ``kills'' the ideal $I_n$ mentioned above, in the ring 
$M{\mathbb R}_{\star}$. 

The contribution of this paper is in calculating the coefficient ring of
$BP{\mathbb R}\langle n\rangle$. This is a non-trivial matter, since 
$I_n$ is certainly not a regular ideal in $M{\mathbb R}_{\star}$. In fact,
it is highly surprising that the spectrum $BP{\mathbb R}\langle n\rangle$ 
constructed in this ``naive'' way gives the coefficients that S. Wilson 
asked for.

To explain the issues involved, 
it should be mentioned at this point that we are dealing here with 
$RO({\mathbb Z}/2)$-graded ${\mathbb Z}/2$-equivariant spectra~\cite{araki1},
\cite{lms}
($M{\mathbb R}$ is ${\mathbb Z}/2$-equivariant), and that, therefore, 
questions of a ``completion theorem''~(\cite{gm}) arise. Indeed, Steve 
Wilson originally asked if a ``homotopy fixed point spectrum'' of 
$BP\langle n\rangle$ is the answer to his question. In this paper,
we shall see that that is, in fact, false. The homotopy fixed point spectrum 
of $BP\langle n\rangle$ will be relevant to our calculations, but
turns out not to have the right coefficients (they contain some 
spurious elements); the point is
that the spectrum $BP{\mathbb R}\langle n
\rangle$ constructed by killing the idea $I_n$ in $M{\mathbb R}$ does not 
satisfy a ``completion theorem'' in the sense of~\cite{gm}.

This also makes our calculation new technically. In~\cite{hk}, where 
coefficients of numerous spectra obtained from $M{\mathbb R}$ by killing 
ideas are calculated, completion theorems for the relevant spectra
always hold and are the bases of
all the calculations. The present paper contains the 
first case where a calculation of coefficient of a ``derived spectrum 
of $M{\mathbb R}$'' is given where the spectrum does not satisfy a completion 
theorem (with the exception of ${\mathbb Z}/2$-equivariant constant Mackey 
functor spectra $H{\mathbb Z}/2$ and $H{\mathbb Z}$, which, in fact, could 
be called $BP{\mathbb R}\langle -1\rangle$ and $BP{\mathbb R}\langle 0
\rangle$ from the point of view of this paper). 
I get my calculations by computing all the 
other terms of the ``Tate diagram'' of Greenlees-May~\cite{gm}. It is somewhat
amazing that the coefficients of $BP{\mathbb R}\langle n\rangle$ defined 
and calculated in this way are a quotient of $M{\mathbb R}_{\star}$, while
the coefficients of 
the other terms of the Tate diagram, notably the Borel cohomology 
spectrum, are not.

In Section 2, I give a short review of Real cobordism theory and the 
main tools used in the paper, as well as the result of the 
calculations for $BP{\mathbb R}\langle n\rangle_{\star}$. In particular, for
$n=1$, we get that the fixed points spectrum $BP{\mathbb 
R}\langle 1\rangle^{{\mathbb Z}/2}$ is $kO$, the connective cover of 
orthogonal $K$-theory $KO$. However, this is not true in other 
twists, i.e., if we first suspend $BP{\mathbb R}\langle 1\rangle$ by copies of 
the sign representation of ${\mathbb Z}/2$ and then take its 
fixed points. In Section 3, I compute the coefficients of
the Tate and Borel cohomology 
spectra of $BP{\mathbb R}\langle n\rangle$, which appear in the 
Tate diagram for 
$BP{\mathbb R}\langle n\rangle $. Finally, 
in Section 4, I calculate the geometric
spectrum of $BP{\mathbb R}\langle n\rangle$ to fill in 
the Tate diagram, and use
it to get the coefficients of $BP{\mathbb R}\langle n\rangle$ itself. It is 
interesting to note that the coefficients of the Tate and geometric
spectra of $BP{\mathbb R}\langle n\rangle$ 
are small: in this sense, one might say that
$BP{\mathbb R}\langle n\rangle$ ``nearly has descent''.

\section{Review of $M{\mathbb R}$-theory and statement of the main result}

I shall 
now describe some basic aspects of Landweber cobordism theory~(\cite{land},
\cite{araki2}, \cite{araki3},
\cite{hk}). First, the infinite loop spaces making up 
$M{\mathbb R}$ are the same as the infinite loop spaces of the 
complex cobordism spectrum $MU$, but there is a ${\mathbb Z}/2$-action,
and the dimensions are indexed 
differently.  Namely,
on the level of prespectra, $MU$ is obtained from the sequence of 
Thom spaces of $n$-dimensional canonical complex bundles $\gamma_n$ on $BU(n)$.
Denoting the Thom space of $\gamma_n$ by 
$BU(n)^{\gamma_n}$, we get structure maps 
\[ \Sigma^{2} BU(n)^{\gamma_n}
\rightarrow BU(n+1)^{\gamma_{n+1}} . \]
In the case of $M{\mathbb R}$, we use the same Thom spaces 
$BU(n)^{\gamma_n}$, but with the ${\mathbb Z}/2$-action by complex 
conjugation.  The space $BU(n)^{\gamma_n}$ is placed in dimension 
$n(1+\alpha)$, where $1$ and $\alpha$ denote the trivial and the 
sign representations of ${\mathbb Z}/2$, respectively. This is because 
$\gamma_n$ is a Real bundle in the sense of Atiyah~\cite{atiyah}. Hence,
the structure maps are 
\[ \Sigma^{1+\alpha}BU(n)^{\gamma_n} \rightarrow BU(n+1)^{\gamma_{n+1}} . \]
%In the case of $M{\mathbb R}$, we also consider the ${\mathbb Z}/2$-action 
%via complex conjugation, and the canonical Real bundles in the 
%sense of Atiyah~\cite{atiyah}.  In particular, let $\alpha$ be 
%the sign representation of ${\mathbb Z}/2$. Then the structure map
%goes from the $(1+\alpha)$-th suspension of the 
%$n$-th space of the prespectrum for $M{\mathbb R}$ to 
%the $(n+1)$-st space. 
(Note that the ${\mathbb Z}/2$-representation 
$1+ \alpha$ is just ${\mathbb C}$ with ${\mathbb Z}/2$-action 
by complex conjugation.)
Spectrification then 
makes $M{\mathbb R}$ a ${\mathbb Z}/2$-equivariant spectrum, 
indexed on $RO({\mathbb Z}/2)$, i.~e. in  dimensions $k + l\alpha$, 
for all $k, l \in {\mathbb Z}$.  In this paper, we will denote 
the $RO({\mathbb Z}/2)$-grading by the subscript $\star$,
to distinguish it from ${\mathbb Z}$-grading, which will be denoted by 
the subscript $\ast$ as usual.  We work locally at the prime 2 in this paper.
The Real Brown-Peterson spectrum is 
obtained from $M{\mathbb R}$ via the Real version of the Quillen idempotent,
analogous to the way the Brown-Peterson spectrum $BP$ is obtained 
from $MU$.  In~\cite{hk}, we calculated the 
$RO({\mathbb Z}/2)$-graded coefficient ring of $BP{\mathbb R}$.
Namely, we have that 
\begin{equation}
BP{\mathbb R}_{\star} = {\mathbb Z}_{(2)}[v_n\sigma^{l2^{n+1}}, a]/
\sim .  
\label{bpr}
\end{equation}
The relations are 
\begin{align}
v_0 & = 2 \\
(v_n\sigma^{l2^{n+1}})a^{2^{n+1}-1} & = 0 \\
(v_n\sigma^{l2^{n+1}})(v_m\sigma^{k2^{m+1}}) & = v_n v_m 
\sigma^{l2^{n+1} + k2^{m+1}} \ \mathrm{for\ } m \leq n .
\end{align}
Here, $n \geq 0$, and $l$ ranges over all integers.
The dimensions of elements are that $v_n$ has dimension $(2^n-1)(1+ \alpha)$,
$a$ has dimension $-\alpha$, and the operator $\sigma$ has dimension 
$-1+\alpha$.  

As described in~\cite{hk},
for each $n \geq 0$, the Real Johnson-Wilson spectrum $BP{\mathbb R}
\langle n\rangle$ is 
obtained by killing the sequence of elements $v_{n+1}, v_{n+2}, \ldots$ in
$BP{\mathbb R}$, in the manner of~\cite{ekmm}. This is again a
${\mathbb Z}/2$-equivariant spectrum indexed on $RO({\mathbb Z}/2)$.
In particular, the 
infinite loop space of $BP{\mathbb R}\langle n \rangle$ in dimension 
$k + l\alpha$ 
is the same as the infinite loop space of $BP\langle n\rangle$ in dimension 
$k + l$, but with an additional action by ${\mathbb Z}/2$, which 
depends on $k$ and $l$, not just their sum.

%The main tool that we will use to compute $BP{\mathbb R}\langle n
%\rangle_{\star}$ 
%is the Tate diagram.  

For a ${\mathbb Z}/2$-equivariant spectrum $E$,
there are several kinds of ``fixed points spectra'' associated with 
$E$.  What we usually consider as the fixed point spectrum is
the Lewis-May fixed point spectrum
$E^{{\mathbb Z}/2}$, obtained by first forgetting 
the $RO({\mathbb Z}/2)$-graded spectrum to one graded on ${\mathbb Z}$, 
i.e., considering only the spaces in dimensions $k + 0 \alpha$, and then 
taking the fixed points spacewise~\cite{lms}. This gives a nonequivariant spectrum.
Similarly, for each $l \in {\mathbb Z}$, one also has 
$(\Sigma^{-l\alpha}E)^{{\mathbb Z}/2}$, called the fixed point 
spectrum twisted by $l$.  This is obtained by first taking only 
the ${\mathbb Z}$-graded spectrum consisting of the spaces in dimensions 
$k + l\alpha$, and then taking fixed points spacewise.  There are also 
the Borel homology and cohomology fixed point spectra of $E$.
Recall that $E{\mathbb Z}/2$ is
the universal contractible free 
${\mathbb Z}/2$-space, which may be thought of as $S(\infty \alpha) =
\colim_k S(k\alpha)$, where $S(k\alpha)$ is the unit sphere in the 
representation $k\alpha$.  The Borel homology spectrum 
$E{\mathbb Z}/2_{+} \wedge E$, and the Borel cohomology 
spectrum is $F(E{\mathbb Z}/2_{+}, E)$.  The Borel 
homology and cohomology fixed points of $E$ are obtained 
by taking the fixed points (in the above sense, with possible twist by 
$l$) of the Borel homology and cohomology spectra, respectively.
In particular, the Borel cohomology fixed points 
\[ F(E{\mathbb Z}/2_{+}, E)^{{\mathbb Z}/2} \]
is $E^{h{\mathbb Z}/2}$, the homotopy fixed points spectrum 
of $E$.  For the Borel homology, note that since 
$E{\mathbb Z}/2_{+} \wedge E$ is a free spectrum indexed 
on $RO({\mathbb Z}/2)$, its fixed points can be computed using the 
Adams isomorphism, which gives that 
\[ (E{\mathbb Z}/2_{+} \wedge E)^{{\mathbb Z}/2} 
\simeq E{\mathbb Z}/2_{+} \wedge_{{\mathbb Z}/2} E . \]

We also have the geometric fixed points spectrum of $E$.
This is is a ${\mathbb Z}$-graded nonequivariant spectrum, whose infinite 
loopspace is 
\[ \colim_{V} \Omega^{V^{{\mathbb Z}/2}}
E_V^{{\mathbb Z}/2}  \]
where the colimit ranges over all finite-dimensional representations $V$
of ${\mathbb Z}/2$, and $E_V$ denotes the $V$-th space 
of $E$.  The geometric fixed points can be calculated 
by first taking $S^{\infty \alpha} \wedge E$, then 
taking the fixed points of this spectrum in the sense above. 
Here, $S^{\infty \alpha}$ 
is the one-point compactification of the infinite-dimensional representaiton 
$\infty \alpha$.

The various spectra associated with a ${\mathbb Z}/2$-equivariant 
spectrum $E$ are organized by 
the Tate diagram.
We have the cofiber sequence 
\[ E{\mathbb Z}/2_{+} \rightarrow S^0 \rightarrow \widetilde{E{\mathbb Z}/2} \]
where the cofiber is the unreduced suspension of $E{\mathbb Z}/2$. Hence, we have that $\widetilde{E{\mathbb Z}/2}$ is just $S^{\infty \alpha}$. 
Smashing with $E$ and mapping into
$F(E{\mathbb Z}/2_{+}, E)$ gives the Tate diagram 
\[ \diagram 
E{\mathbb Z}/2_{+} \wedge E \dto_{\simeq} \rto & 
E \dto \rto & \widetilde{E{\mathbb Z}/2} \wedge 
E \dto \\
E{\mathbb Z}/2_{+} \wedge F(E{\mathbb Z}/2_{+}, E) \rto &
F(E{\mathbb Z}/2_{+}, E) \rto & \widetilde{E{\mathbb Z}/2}
\wedge F(E{\mathbb Z}/2_{+}, E) .
\enddiagram \]
The rightmost term on the bottom row, $\widetilde{E{\mathbb Z}/2} \wedge
F(E{\mathbb Z}/2_{+}, E)$, is the Tate cohomology of 
$E$, which we also denote by $t(E)$.

Taking $E = BP{\mathbb R}\langle n\rangle$, we get the Tate diagram for 
$BP{\mathbb R}\langle n \rangle$
\[ \diagram 
E{\mathbb Z}/2_{+} \wedge BP{\mathbb R}\langle n\rangle \dto_{\simeq} \rto & 
BP{\mathbb R}\langle n\rangle \dto \rto & \widetilde{E{\mathbb Z}/2} \wedge 
BP{\mathbb R}\langle n\rangle \dto \\
{\begin{array}{c} E{\mathbb Z}/2_{+}  \\
\wedge F(E{\mathbb Z}/2_{+}, BP{\mathbb R}\langle 
n\rangle) 
\end{array}} \rto &
F(E{\mathbb Z}/2_{+}, BP{\mathbb R}\langle n\rangle) \rto & 
t(BP{\mathbb R}\langle n\rangle) .
\enddiagram \]
Here, $t(BP{\mathbb R}\langle n\rangle) =
\widetilde{E{\mathbb Z}/2}
\wedge F(E{\mathbb Z}/2_{+}, BP{\mathbb R}\langle n\rangle))$. 
One sometimes also refers to the fixed points spectra obtained 
from the spectra in the Tate diagram by the same names as the
corresponding equivariant spectra.  Note that we can also take 
twisted fixed points, by first desuspending by $S^{l\alpha}$, and 
then taking fixed points.  However, note that the rightmost column,
i.e., the geometric and the Tate spectra, are $\alpha$-periodic, and 
hence do not depend on the twist $l$. The middle column of the Tate 
diagram is an equivalence if and only if the rightmost column is an 
equivalence.  We call an $RO({\mathbb Z}/2)$-graded
equivariant spectrum \emph{complete} if this condition holds. More 
generally, a ``completion theorem'' holds if the middle vertical arrow 
of the Tate diagram is a completion in some suitable sense. For 
more information, see~\cite{gm}.  
Unlike $BP{\mathbb R}$, the spectrum $BP{\mathbb R}\langle n\rangle$ is 
not complete, i.e., it is not equivalent to its Borel cohomology 
spectrum. 

All these spectra help in computing the coefficients of 
$BP{\mathbb R}\langle n\rangle$. There are spectral sequences that compute the
coefficients of the Borel homology, Borel cohomology, and Tate cohomology 
terms, while $\widetilde{E{\mathbb Z}/2} \wedge BP{\mathbb R}
\langle n\rangle$ can 
be identified using geometric methods.
%The three terms of the bottom row are the Borel homology spectrum, 
%the Borel cohomology spectrum, and the Tate spectrum of $BP{\mathbb R}
%\langle n\rangle$,
%respectively.

\begin{theor}
{\rm (1)}\qua  The coefficients of the Tate spectrum of $BP{\mathbb R}\langle n
\rangle$ are
\[ t(BP{\mathbb R}\langle n \rangle)_{\star} = {\mathbb Z}/2[\sigma^{2^{n+1}}, \sigma^{-2^{n+1}}, a, a^{-1}] \]
where $\sigma$ has dimension $-1 + \alpha$, and $a$ has dimension $-\alpha$.

{\rm (2)}\qua The coefficients of the Borel cohomology spectrum of $BP{\mathbb R}
\langle n\rangle$ are 
\[ F(E{\mathbb Z}/2_{+}, BP{\mathbb R}\langle n \rangle)_{\star} 
= ({\mathbb Z}_{(2)}[v_k \sigma^{l2^{k+1}}, a]/\sim)
\oplus {\mathbb Z}/2[\sigma^{2^{n+1}}, \sigma^{-2^{n+1}}, a] . \]
Here, $0 \leq k \leq n$, and $l$ ranges over all integers. The relations 
are 
\begin{align*}
v_0 & = 2 \\
v_k a^{2^{k+1}-1} & = 0 \\
(v_n\sigma^{l2^{n+1}})(v_m\sigma^{k2^{m+1}}) & = 
v_nv_m\sigma^{l2^{n+1} + k2^{m+1}} \ \mathrm{for\ } m \leq n .
\end{align*}
\label{tateborel}
\end{theor}

For $BP{\mathbb R}\langle n\rangle$ itself, we have the following theorem.
\begin{theor}
The coefficient ring of $BP{\mathbb R}\langle n\rangle$ is 
\[ BP{\mathbb R}\langle n \rangle_{\star} = ({\mathbb Z}_{(2)}[v_k \sigma^{l2^{k+1}}, a]/\sim)
\oplus {\mathbb Z}/2[\sigma^{-2^{n+1}}, a] . \]
\label{bprn1}
with the same relations (2.2), (2.3) and (2.4) as in $BP{\mathbb R}_{\star}$.
\end{theor}

For readers who prefer not to use the $RO({\mathbb Z}/2)$-grading, 
the (untwisted or twisted) coefficients of $BP{\mathbb R}$ and 
$BP{\mathbb R}\langle n\rangle$ can be described using nonequivariant Milnor 
words with dimensional shifts.
For an element $x$ of dimension $k + l \alpha$, we say that the 
twist of $x$ is $l$. Recalling the calculation of $BP{\mathbb R}_{\star}$, 
for a fixed twist $l$, we can describe the coefficients of 
$(\Sigma^{-l\alpha}BP{\mathbb R})^{{\mathbb Z}/2}$, the 
twist $l$ fixed points of $BP{\mathbb R}$, in terms of just the 
Milnor generators $v_n$'s, but with shifted dimensions. Namely, 
fix $l \in {\mathbb Z}$.  For a sequence of nonnegative 
integers $R = (r_0, r_1, \ldots)$ of which all but finitely 
many are $0$, we write the monomial 
$v_R = \prod_{i \geq 0} v_i^{r_i}$. Let $n = min(R)$ be the smallest 
number such that $i_n \neq 0$, and let $|v_R|$ 
denote the dimension of $v_R$ in $BP_{\ast}$. 
The additive generators of $BP{\mathbb R}_{\star}$ as a 
${\mathbb Z}_{(2)}$-module are the monomials $v_R$, with the 
following possibilities.  If $|v_R| \leq l$, the $v_R$ has dimension 
 \[ |v_R| -l - k \]
where $0 \leq k \leq 2^{n+1}-1$ is congruent to $|v_R|$ modulo 
$2^{n+1}$. This is $0$ if $k = 2^{n+1}-1$,
it generates a copy of ${\mathbb Z}_{(2)}$ if $k = 0$,
and a copy of ${\mathbb Z}/2$ otherwise.
%relation that $v_R = 0$ if $|v_R|-l$ \equiv -1$ modulo $2^{n+1}$.
%If $|v_R| -l \not\equiv -1$ modulo 
%$2^{n+1}$, 
%the 
%$v_R = 0$ in $(\Sigma^{-l\alpha}BP{\mathbb R})^{{\mathbb Z}/2}_{\ast}$.
%Otherwise, 
%then we have two possibilities in 
%$(\Sigma^{-l\alpha}BP{\mathbb R})^{{\mathbb Z}/2}_{\ast}$.
If $l > |v_R|$, then 
$v_R$ is in dimension 
\[ |v_R| -l -k^{\prime} \]
where $0 \leq k^{\prime} \leq 2^{n+1}-1$ is congruent to $|v_R| -l$ 
modulo $2^{n+1}$.  Again, this is $0$ if $k^{\prime} = 2^{n+1}-1$,
it generates a copy of ${\mathbb Z}_{(2)}$
if $k^{\prime} = 0$, and it generates a copy of ${\mathbb Z}/2$ else.

For each $l$, the elements of the homotopy groups 
of the twist $l$ fixed points of 
$BP{\mathbb R}\langle n\rangle$ are the 
relevant ones from $BP{\mathbb R}$, and some 
extra elements.
\begin{coro}
Let $l \in {\mathbb Z}$.  If $l \geq 0$, then the elements of 
$BP{\mathbb R}\langle n\rangle_{\star}$ in twist $l$ are the same as the 
twist $l$ 
elements of $BP{\mathbb R}_{\star}$ that do not contain $v_s$ for 
any $s > n$.  If $l<0$, then the elements of 
$BP{\mathbb R}\langle n\rangle_{\star}$ in twist $l$ are the twist $l$ 
elements of
$BP{\mathbb R}_{\star}$ not containing $v_s$ for any $s > n$, 
as well as an extra copy of ${\mathbb Z}/2$ in 
dimension $k2^{n+1}$ for each $k$ such that $ 0 > k2^{n+1} \geq l$.
(In the notation of Theorem 2, this element corresponds to the  
generator by $\sigma^{-k2^{n+1}}a^{k2^{n+1}-l}$.)
\end{coro}

\section{Tate and Borel cohomology calculations}

The goal of this section is to prove Theorem~\ref{tateborel}.  To compute 
the Tate cohomology of $BP{\mathbb R}\langle n\rangle$, we consider the 
Tate spectral sequence 
\begin{equation}
E_2 = \widehat{H}^{\ast}({\mathbb Z}/2, BP\langle n\rangle_{\ast}[\sigma, 
\sigma^{-1}])
\Rightarrow \widehat{BP{\mathbb R}\langle n\rangle}_{\star} . 
\label{bprntate}
\end{equation}
We can compare this to the Tate spectral sequence for $BP{\mathbb R}$, 
which is that 
\begin{equation}
E_2 = \widehat{H}^{\ast}({\mathbb Z}/2, BP_{\ast}[\sigma, \sigma^{-1}])
\Rightarrow \widehat{BP{\mathbb R}}_{\star} . 
\label{bprtate}
\end{equation}
(see~\cite{hk, ext0}).  The $E_1$-term of~(\ref{bprtate}) is 
$BP_{\star}[\sigma, \sigma^{-1}, a, a^{-1}]$, where $BP_{\star}$ is 
the same as $BP_{\ast}$, with the exception that the dimension of 
$v_n$ is $(2^n-1)(1+\alpha)$ instead of $2(2^n-1)$.  
Note that with a different choice of generators (multiplying 
$v_n$ by $\sigma^{2^n-1}$), this is in fact 
equal to $BP_{\ast}[\sigma, \sigma^{-1}, a, a^{-1}]$. In~(\ref{bprtate}),
${\mathbb Z}/2$ acts by $(-1)^{(|v_R|_{\mathbb C})/2 + l}$ on the 
monomial $v_R\sigma^{l}$.  Here, $|v_R|_{\mathbb C}$ denotes the 
dimension of a monomial $v_R$ in $BP_{\ast}$.
In~\cite{hk}, it was shown that~(\ref{bprtate})
has the differentials 
\begin{equation}
 d_{2^{k+1}-1}(\sigma^{-2^k}) = v_ka^{2^{k+1}-1} 
\label{differential}
\end{equation}
for $k \geq 1$.
These differentials wipe out all elements except ${\mathbb Z}/2[a, a^{-1}]$.
Namely, a typical element of the $E_1$-term of the spectral 
sequence~(\ref{bprtate}) is $v_R\sigma^{2^sl}a^t$, where $l \in 
{\mathbb Z}$ is odd, $t \in {\mathbb Z}$, and $R = (r_0, r_1, \ldots)$
is a sequence of nonnegative integers, of which only finitely many are 
nonzero, with $v_R = \prod v_i^{r_i}$.  The filtration degree of this 
element is $t$.  The differential~(\ref{differential})
gives that if $s \leq min(R)$, then 
\begin{equation}
 d_{2^{s+1}-1}(v_R \sigma^{2^sl}a^t) = v_R v_s \sigma^{2^s{l+1}}a^{t+
2^{s+1}-1}  
\label{tatediff}
\end{equation}
for all $l \neq 0$.
So the element is the source of a differential if $s \leq min(R)$ or 
if $R = (0, 0, \ldots)$ and $l \neq 0$, and it is the target of a 
differential if $s > min(R)$ or if $l = 0$ and $R \neq (0, 0, \ldots)$.
Note that every monomial $v_{R}\sigma^{2^sl}a^t$ in the $E_1$-term of 
appears either in the source or target of a differential~(\ref{tatediff}),
except when $l = 0$. (For complete details on this, see~\cite{ext0}.)
Thus, the only surviving elements are powers of $a$.

In~(\ref{bprntate}), the $E_1$-term is now 
\[ BP\langle n\rangle_{\star}[\sigma, \sigma^{-1}, a, a^{-1}] . \]
Again, this is the same as $BP\langle n\rangle_{\ast}[\sigma, 
\sigma^{-1}, a, a^{-1}]$, by replacing the generators $v_i$ by 
$v_i\sigma^{2^i-1}$. 
The differentials are same as the ones as~(\ref{differential}). 
An element of the $E_1$-term is $v_R\sigma^{2^sl}a^t$, but now 
$R = (v_0, v_1, \ldots, v_n)$. If $s \leq min(R)$, this is the 
source of a differential. If $s > min(R)$ or if $l = 0$ and $R \neq 
(0, 0, \ldots, 0)$, this is the 
target of a differential.  However, suppose that $R = (0, \ldots, 0)$ and 
$s > n$. In the spectral sequence~(\ref{bprntate}), we get a differential
\begin{equation}
 d_{2^{s+1}-1}(\sigma^{2^sl}a^t) = v_s \sigma^{2^s(l+1)}a^{t+
2^{s+1}-1} . 
\end{equation}
The target of this differential is now $0$ in the Tate spectral sequence 
for $B\!P{\mathbb R}\langle n\rangle$. 
%For $0 \leq k \leq n$, these differentials exist, and wipe out all
%monomials with $v_k$'s, as well as all $\sigma^j$ such that 
%$j$ is not a multiple of $2^{n+1}$.  
%In~(\ref{bprtate}), the element $\sigma^{-2^{n+1}}$ survives to 
%$E_{2^{n+1}-1}$, and is then a source of $d_{2^{n+1}-1}$. In this 
%case, however, the differential 
%since $v_{n+1}$ is not contained in $BP\langle n\rangle_{\ast}$. 
Thus, the monomials in
$\sigma^{2^{n+1}}, \sigma^{-2^{n+1}}, a, a^{-1}$ survive to the 
$E_{\infty}$-term of the Tate spectral sequence~(\ref{bprntate}). This 
proves the first part of Theorem~\ref{tateborel}.

For the Borel cohomology of $BP{\mathbb R}\langle n\rangle$, we use the 
Borel cohomology spectral sequence 
\begin{equation}
E_2 = H^{\ast}({\mathbb Z}/2, BP\langle n\rangle_{\ast}[\sigma, 
\sigma^{-1}]) \Rightarrow 
F(E{\mathbb Z}/2_{+}, BP{\mathbb R}\langle n\rangle)_{\star} . 
\label{bprnborel}
\end{equation}
We compare this to both the Tate spectral sequence~(\ref{bprntate}), 
and to the Borel cohomology spectral sequence for $BP{\mathbb R}$
\begin{equation}
E_2 = H^{\ast}({\mathbb Z}/2, BP_{\ast}[\sigma, \sigma^{-1}]) \Rightarrow 
F(E{\mathbb Z}/2_{+}, BP{\mathbb R})_{\star} . 
\label{bprborel}
\end{equation}
The $E_1$-term of~(\ref{bprnborel}) is $BP\langle n\rangle_{\star}[\sigma, 
\sigma^{-1}, a]$,
which is just the part of the $E_1$-term of the Tate spectral 
sequence~(\ref{bprntate}) consisting of only the elements with nonnegative 
filtration degrees (i.~e. nonnegative powers of $a$).  The differentials are 
the same as in~(\ref{bprntate}), 
i.~e. $d_{2^{k+1}-1}$ for $0 \leq k \leq n$, except that now we only 
allow the differentials with sources and targets both having nonnegative 
filtration degrees.  Thus, the monomials $v_R \sigma^{2^{s}l}a^t$ 
with $t \leq 2^{min(R)+1}-2$ and $s > min(R)$ will survive~(\ref{bprnborel}), 
since in~(\ref{bprntate}), they are targets of differentials $d_{2^{s+1}-1}$
with sources having negative filtration degrees. Also, as before,
the monomials $\sigma^{2^sl} a^t$ survives for any $s > n$ and 
$t \in {\mathbb Z}$.  This gives the second 
part of Theorem~\ref{tateborel}.

\section{The coefficients of $BP{\mathbb R}\langle n\rangle$}

We prove Theorem~\ref{bprn1} in this section. To this end, we will first 
compute the coefficients of the geometric spectrum 
$\widetilde{E{\mathbb Z}/2} \wedge BP{\mathbb R}\langle n\rangle$ by induction 
on $n$. The following lemma was shown in~\cite{hk}.

\begin{lemma}
$\widetilde{E{\mathbb Z}/2} \wedge
BP{\mathbb R}\langle 0\rangle$ is $H{\mathbb Z}/2_{m}$, the 
$\mathbb{Z}/2$-equivariant 
cohomology spectrum corresponding to the constant Mackey functor.
\end{lemma}

\begin{prop}
For $n \geq 0$, the coefficients of the geometric
spectrum $\widetilde{E{\mathbb Z}/2}$ $\wedge\, BP{\mathbb R}\langle 
n\rangle$ are 
\[ {\mathbb Z}/2[\sigma^{-2^{n+1}}, a, a^{-1}] \]
where the dimensions of $\sigma$ and $a$ are as above.
\end{prop}

\begin{proof}
We work by induction.  As shown in~\cite{hk}, the 
coefficents of $\widetilde{E{\mathbb Z}/2} \wedge BP{\mathbb R}
\langle 0\rangle$ is 
\[ (\widetilde{E{\mathbb Z}/2} \wedge 
BP{\mathbb R}\langle 0\rangle)_{\star} = (H{\mathbb Z}/2_m)_{\star} = 
{\mathbb Z}/2[\sigma^{-2}, a, a^{-1}] . \]
Suppose that the statement is true for $\widetilde{E{\mathbb Z}/2} 
\wedge BP{\mathbb R}\langle n-1\rangle$.
We filter $\widetilde{E{\mathbb Z}/2} \wedge 
BP{\mathbb R}\langle n\rangle$ by copies of $\widetilde{E{\mathbb Z}/2} \wedge 
BP{\mathbb R}\langle n-1\rangle$. Namely, consider the map 
\[ v_n: \Sigma^{(2^n-1)(1+\alpha)}(\widetilde{E{\mathbb Z}/2} 
\wedge BP{\mathbb R}\langle n\rangle) \rightarrow 
\widetilde{E{\mathbb Z}/2} \wedge BP{\mathbb R}
\langle n\rangle . \]
The cofiber of this is a suspension of $\widetilde{E{\mathbb Z}/2} 
\wedge BP{\mathbb R}\langle n-1\rangle$.  Iterating the map gives an 
exact couple, which in turn gives a spectral sequence
\begin{equation}
\begin{split}
E_1 = (\widetilde{E{\mathbb Z}/2} \wedge
BP{\mathbb R}\langle n-1\rangle)_{\star}[v_n] & = 
{\mathbb Z}/2[\sigma^{-2^{n}}, 
a, a^{-1}][v_n] \\
& \Rightarrow (\widetilde{E{\mathbb Z}/2} 
\wedge BP{\mathbb R}\langle n\rangle)_{\star}. 
\end{split}
\label{geoss}
\end{equation}
By comparing with the other spectral sequences~(\ref{bprntate}),
(\ref{bprtate}) and~(\ref{bprnborel}), we get that the differentials 
of the spectral sequence are only 
\[ d_1(\sigma^{-2^n}) = v_n a^{2^{n+1}-1} \]
and its multiples by powers of $\sigma^{-2^n}$.   
Hence, by arguments similar to that for the Tate spectral sequence,
$\sigma^{-l2^n}$ is the source of a differential for all $l$ odd, and 
all monomials containing $v_n$ are targets of differentials.
This gives 
that the $E_{\infty}$-term of the spectral sequence~(\ref{geoss}) is 
just ${\mathbb Z}/2[\sigma^{-2^{n+1}}, a, a^{-1}]$ as claimed.
\end{proof}

Now in the bottom row of the Tate diagram for $BP{\mathbb R}
\langle n\rangle$, 
we have the cofiber sequence 
\begin{equation*}
\begin{split}
& E{\mathbb Z}/2_{+} \wedge F(E{\mathbb Z}/2_{+}, BP{\mathbb R}
\langle n\rangle)_{\star} \\
& \rightarrow ({\mathbb Z}_{(2)}[v_k\sigma^{l2^{n+1}}, a]/\sim) \oplus
{\mathbb Z}/2 [\sigma^{2^{n+1}}, \sigma^{-2^{n+1}}, 
a] \\
& \rightarrow {\mathbb Z}/2[\sigma^{2^{n+1}}, 
\sigma^{-2^{n+1}}, a, a^{-1}] .
\end{split}
\end{equation*}  
By comparison of the spectral sequences 
computing them, it is straightforward to see that the map from the 
Borel cohomology term to the Tate term is just the inclusion on the 
monomials containing only $a$ and powers of $\sigma$, and kills all 
monomials containing any $v_k$. 
Thus, the coefficient of the fibers is
\[ ({\mathbb Z}_{(2)}[v_k\sigma^{l2^{n+1}}, a]/\sim) \oplus {\mathbb Z}/2
[a^{-1}] . \]
Here, $\sim$ denotes the relations (2.2), (2.3) and (2.4).
This is the Borel homology of $BP{\mathbb R}
\langle n\rangle$. Hence, for the top row 
of the Tate diagram, we get the cofiber sequence
\begin{equation*}
\begin{split}
({\mathbb Z}_{(2)}[v_k\sigma^{l2^{n+1}}, a]/\sim) \oplus {\mathbb Z}/2
[a^{-1}] & \rightarrow BP{\mathbb R}\langle n\rangle_{\star} \\
& \rightarrow {\mathbb Z}/2[\sigma^{-2^{n+1}}, a, a^{-1}] .
\end{split}
\end{equation*}
The connecting map is the identity on $a^{-1}$ and kills $\sigma^{-2^{n+1}}$
and $a$. Therefore, the middle term gives 
\[ BP{\mathbb R}\langle n\rangle_{\star} 
= ({\mathbb Z}_{(2)}[v_k\sigma^{l2^{n+1}}, a]/\sim)
\oplus {\mathbb Z}/2[\sigma^{-2^{n+1}}, a]  \]
where $\sim$ denotes the relations (2.2), (2.3) and (2.4).

\np

\Addresses\recd

\end{document}